\documentstyle[11pt]{article}
\input epsf.tex

\newtheorem{example}{Example}[section]

\newtheorem{theorem}[example]{Theorem}

\newtheorem{proposition}[example]{Proposition}

\newtheorem{conjecture}[example]{Conjecture}

\def\boxit#1#2{\setbox1=\hbox{\kern#1{#2}\kern#1}%
\dimen1=\ht1 \advance\dimen1 by #1 \dimen2=\dp1 \advance\dimen2 by #1
\setbox1=\hbox{\vrule height\dimen1 depth\dimen2\box1\vrule}%
\setbox1=\vbox{\hrule\box1\hrule}%
\advance\dimen1 by .4pt \ht1=\dimen1
\advance\dimen2 by .4pt \dp1=\dimen2 \box1\relax}

\def\adots{\mathinner{\mkern2mu\raise1pt\hbox{.}
\mkern3mu\raise4pt\hbox{.}\mkern1mu\raise7pt\hbox{.}}}
\def\<{\langle}
\def\>{\rangle}

\font\tensym=msbm10
\font\sevensym=msbm7
\font\fivesym=msbm5

\newfam\ssymfam
\textfont\ssymfam=\tensym
\scriptfont\ssymfam=\sevensym
\scriptscriptfont\ssymfam=\fivesym
\def\ssym{\fam\ssymfam\tensym}

\def\Z{{\ssym Z}}
\def\Q{{\ssym Q}}
\def\C{{\ssym C\,}}

\font\tengoth=eufm10
\font\sevengoth=eufm7
\font\fivegoth=eufm5

\newfam\gothfam
\textfont\gothfam=\tengoth
\scriptfont\gothfam=\sevengoth
\scriptscriptfont\gothfam=\fivegoth
\def\goth{\fam\gothfam\tengoth}

\def\I{{\cal J}}

\newfont{\bb}{cmbx10}

\def\mod{{\rm \ mod\ }}
\def\D{{\bf D}}

\def\U{{\bf U}}

\def\glchap{\widehat{\hbox{\goth gl}}}
\def\slchap{\widehat{\hbox{\goth sl}}}

\def\F{{\cal F}}

\def\H{{\cal H}}
\def\F{{\cal F}}

\def\F{{\cal F}}

\title{\Large\bf Canonical Bases of $q$-Deformed Fock Spaces}

\author{\rm Bernard \sc Leclerc\thanks{Universit\'e de Caen, D\'epartement
de Math\'ematiques, Esplanade de la Paix, BP 5186, 14032 Caen cedex, France}
and
\rm Jean-Yves \sc Thibon\thanks{Institut Gaspard Monge, Universit\'e de
Marne-la-Vall\'ee, 2 rue de la Butte-Verte, 93166 Noisy-le-Grand cedex, France}}
\date{}

\begin{document}

\maketitle

\begin{abstract}
We define a canonical basis of the $q$-deformed Fock space representation
of the affine Lie algebra $\glchap_n$. We conjecture that
the entries of the transition matrix between this basis and the natural 
basis of the Fock space are $q$-analogues of decomposition numbers
of the $v$-Schur algebras for $v$ specialized
to a $n$th root of unity.
\end{abstract}

%%%%%%%%%%%%%%%%%%%%%%%%%%%%%%%%%%%%%%%%%%%%%%%%%%%%%%%%%%%%%%%%%%%%%%%%%%%

\section{Introduction}
The Fock space representation ${\cal F}$ of $\slchap_n$
is not irreducible. Its decomposition into simple $\slchap_n$-modules
is given by \cite{DJKM}
\begin{equation}\label{DECOMP}
\F \cong \bigoplus_{k\ge 0} M(\Lambda_0 - k\delta)^{\oplus p(k)} \,,
\end{equation}
where $p(k)$ denotes the number of partitions $\lambda$ of $k$.
Hence it is not obvious to apply 
Kashiwara's or Lusztig's method to define a canonical basis
of ${\cal F}$.

In this note, we shall rather regard ${\cal F}$ as a representation
of the enlarged algebra $\glchap_n$. Indeed, 
${\cal F}$ is a simple $\glchap_n$-module, and a $q$-deformation
${\cal F}_q$ of this representation has been described in \cite{KMS}.
One can then define a natural semi-linear
involution $v \rightarrow \overline{v}$ commuting with
the action of the lowering operators of ${\cal F}_q$ and leaving invariant 
its highest weight vector. 
Using this involution, one obtains in an elementary way
a canonical basis $\{G(\lambda)\}$ of ${\cal F}_q$. 
This basis can be computed explicitely.
It would be interesting to compare it with the canonical basis
obtained via a geometric approach by Ginzburg, Reshetikhin and 
Vasserot~\cite{GRV}.

By restriction to $U_q(\slchap_n)$, the space ${\cal F}_q$
decomposes similarly as
\begin{equation}\label{QDECOMP}
\F_q \cong \bigoplus_{k\ge 0} M_q(\Lambda_0 - k\delta)^{\oplus p(k)} \,,
\end{equation} 
where $M_q(\Lambda)$ denotes the simple $U_q(\slchap_n)$-module
with highest weight $\Lambda$.
The $G(\mu)$ indexed by $n$-regular partitions $\mu$
coincide with the elements of Kashiwara's global crystal basis of 
the basic representation $M(\Lambda_0)$.
But we insist that the rest of our basis is not compatible with the
decomposition (\ref{QDECOMP}). 

We conjecture that for $q=1$, the coefficients of the transition matrix
of our canonical basis on the natural basis of the Fock space
are equal to the decomposition numbers of $v$-Schur algebras
over a field of characteristic 0
at a $n$-th root of unity. A previous conjecture \cite{LLT1,LLT2}
on decomposition matrices of Hecke algebras having been recently
established by Ariki and by Grojnowski, we already know that
the columns of the transition matrix indexed by $n$-regular
partitions contain $q$-analogues of decomposition numbers.
On the other hand, we can prove by means
of a $q$-analogue of Steinberg's tensor product theorem that
an infinite number of entries of the inverse transition matrix
are $q$-analogues of inverse decomposition numbers.

Details and proofs will appear in a forthcoming paper.

We follow the notation of \cite{Mcd} for symmetric functions,
and that of \cite{St,KMS} for $q$-wedges and Fock space representations,
except for the replacement of $q$ by $q^{-1}$.

%%%%%%%%%%%%%%%%%%%%%%%%%%%%%%%%%%%%%%%%%%%%%%%%%%%%%%%%%%%%%%%%%%%%%%%%%%%%
%%%%%%%%%%%%%%%%%%%%%%%%%%%%%%%%%%%%%%%%%%%%%%%%%%%%%%%%%%%%%%%%%%%%%%%%%%%%

\section{A $q$-analogue of the Fock space representation of $\glchap_n$}

The Lie algebra $\glchap_n$ can be regarded  as the sum $\slchap_n+\H_n$
where $\H_n$ is a Heisenberg algebra  commuting with $\slchap_n'$ and
such that $\slchap_n'\cap\H_n=\C c$, where $c$ is the central generator
\cite{DJKMO}. 

The bosonic Fock space $\F=\C[x_1,x_2,\ldots]$ can be interpreted
as the algebra of symmetric functions  in some infinite set of variables
$t_i$ by means of the correspondence $x_k={1\over k}p_k$, where
$p_k=\sum_i t_i^k$ are the power sums. 
With this interpretation, the action of $\glchap_n$ on $\F$ can be 
described as follows. The generator $B_k$ of $\H_n$ acts by
$nk{\partial \over \partial p_{nk}}$ for $k>0$ and as the multiplication
by $p_{-nk}$ for $k<0$. The action of the generators of $\slchap_n$
is particularly simple in the basis of Schur functions $s_\lambda$.
For a node $\gamma$ of the Young diagram of a partition $\lambda$,
located at the intersection of the $i$th row and the $j$th column of
$\lambda$, define its residue  
$r(\gamma)\in\{0,1,\ldots,n-1\}$ as $r(\gamma)=j-i \mod n$.
Then, 
$$e_i s_\lambda = \sum s_\nu\,,\quad \quad  f_i s_\lambda = \sum s_\mu\,,$$
where $\nu$ (resp. $\mu$) runs through all partitions obtained from
$\lambda$ by removing (resp. adding) a node of residue $i$.

The $q$-analogue $\F_q$ of the Fock space representation of $\slchap_n$,
in the form of a Fock space representation of the quantized enveloping
algebra $U_q(\slchap_n)$,
has been constructed  by Hayashi \cite{Ha} and further investigated
by Misra and Miwa \cite{MM} who constructed the crystal basis. 
Recently, Kashiwara, Miwa and Stern~\cite{KMS} have shown that
the action of $U_q(\slchap_n')$ on $\F_q$ is centralized by
a Heisenberg algebra $\H_n^q$.
Let $\U$ be the subalgebra of ${\rm End}(\F_q)$ generated by
these actions of $U_q(\slchap_n)$ and $\H_n^q$ (it would be interesting
to compare $\U$ with the standard $q$-deformation $U_q(\glchap_n)$
considered in \cite{DF,GRV}). 
The Fock space is an irreducible $\U$-module, and
the canonical basis constructed in this paper will be adapted
to this representation.

The Fock space representation of $U_q(\slchap_n)$ can be described
as follows (note that our conventions are slightly different from those
of \cite{MM} and \cite{KMS}).  Let 
$$
\F_q=\bigoplus_\lambda \Q(q) |\lambda\>
$$
be the $\Q(q)$ vector space with basis
$|\lambda\>$ indexed by the set of all partitions.
Let $\lambda$ and $\mu$ be two partitions
such that $\mu$ is obtained from $\lambda$ by adding a node $\gamma$ of
residue $i$. Let $I_i(\lambda)$ be the number of indent $i$-nodes
of $\lambda$, $R_i(\lambda)$ the number of its removable $i$-nodes,
$I_i^l(\lambda,\mu)$ (resp. $R_i^l(\lambda,\mu)$)
the number of indent $i$-nodes (resp. of removable $i$-nodes)
situated to the left of $\gamma$ ($\gamma$ not included),
and similarly, let $I_i^r(\lambda,\mu)$ and $R_i^r(\lambda,\mu)$ be  the
corresponding numbers for nodes located on the right of $\gamma$.
Set $N_i(\lambda)=I_i(\lambda)-R_i(\lambda)$,
$N_i^l(\lambda,\mu)=I_i^l(\lambda,\mu)-R_i^l(\lambda,\mu)$ and
$N_i^r(\lambda,\mu)=I_i^r(\lambda,\mu)-R_i^r(\lambda,\mu)$. Then,
\begin{equation}
f_i |\lambda\> = \sum_\mu q^{N_i^r(\lambda,\mu)} |\mu\>\, ,\qquad
e_i |\mu\>     = \sum_\lambda q^{N_i^l(\lambda,\mu)} |\lambda\>
\end{equation}
where in each case the sum is over all partitions such that
$\mu/\lambda$ is a $i$-node,
\begin{equation}
q^{h_i} |\lambda\> = q^{N_i(\lambda)}|\lambda\> \quad
\mbox{\rm and}
\quad q^D |\lambda\> = q^{-N^0(\lambda)} |\lambda\>
\end{equation}
where $D$ is the degree generator and $N^0(\lambda)$ the total
number of $0$-nodes of $\lambda$.

The Fock representation of the Heisenberg algebra $\H_n^q$ of \cite{KMS} is defined by
means of Stern's construction of semi-infinite $q$-wedges \cite{St}.
The basis vector $|\lambda\>$ is interpreted as the semi-infinite
wedge $u_I= u_{i_1}\wedge u_{i_2}\wedge \cdots $ where the sequence
$I$ is defined by $i_k = \lambda_k-k+1$. 
With our conventions, the commutation relations for $q$-wedges
are the following. Suppose that $\ell <m$ and that $\ell-m \mod n = i$.
If $i=0$ then $u_\ell\wedge u_m=-u_m\wedge u_\ell$ and otherwise,
$$
u_\ell\wedge u_m = -q^{-1} u_m \wedge u_\ell
+\left ( q^{-2} -1 \right)
\left(u_{m-i}\wedge u_{\ell+i} - q^{-1} u_{m-n}\wedge u_{\ell+n}
+ q^{-2} u_{m-n-i}\wedge u_{m+n+i} -\cdots \right)
$$
where in the last sum one retains only the normally ordered terms.

Then, the generator $B_k$ of $\H_n^q$ acts on $u_I$ by
\begin{equation}
B_k u_I = \sum_{i\ge 1} u_{I+kn \delta^i}
\end{equation}
where $\delta^i$ is the sequence $(\delta^i_j)_{j\ge 1}$ (Kronecker
symbols). The semi-infinite wedge $u_I$ will be called the fermionic
realization of the basis vector $|\lambda\>$. 
These operators verify the relations \cite{KMS}
\begin{equation}
[B_k,B_\ell] = k {1-q^{-2nk}\over 1-q^{-2k}} \delta_{k+\ell,0} \ .
\end{equation}

The action of $\H_n^q$ in the bosonic picture is better described
in terms of other operators $U_k,V_k\in U(\H_n^q)$ \cite{LLT3}.
Recall that for $k>0$, $B_{-k}$ is a $q$-analogue of multiplication
by $p_{kn} = p_k(t_1^n,t_2^n,\ldots)=\psi^n(p_k)$, where $\psi_n$
is the ring homomorphism raising the variables to the $n$th power.
On the graphical representation by
Young diagrams,  multiplication by the complete homogeneous functions
$\psi_n(h_k)=h_k(t_1^n,t_2^n,\ldots)$ of $n$th powers  has a simple
combinatorial description. 
Let
\begin{equation}
V_k = \sum_{m_1+2m_2+\cdots+km_k=k}
{1\over m_1!m_2!\cdots m_k!}
\left( {B_{-1}\over 1} \right)^{m_1}
\left( {B_{-2}\over 2} \right)^{m_2}\cdots
\left( {B_{-k}\over k} \right)^{m_k}
\end{equation}
be the $q$-analogue of the multiplication operator by $\psi^n(h_k)$.
Then,
\begin{equation}\label{eq8}
V_k |\lambda\> =\sum q^{-{\bf h}(\mu/\lambda)} |\mu\>
\end{equation}
where the sum is over all partitions $\mu$ such that $\mu/\lambda$ is
a horizontal $n$-ribbon strip of weight $k$, and 
$$
{\bf h}(\mu/\lambda)=\sum_R ({\rm ht\,} (R) -1)
$$
where the sum is over all the $n$-ribbons $R$ tiling $\mu/\lambda$,
${\rm ht\,}(R)$ being the height of the ribbon $R$ (see \cite{LLT3}).

The scalar product on $\F_q$ is defined by $\<\lambda|\mu\>=\delta_{\lambda\mu}$.
The adjoint operator $U_k$ of $V_k$ acts by
\begin{equation}
U_k |\mu\> =\sum q^{-{\bf h}(\mu/\lambda)} |\lambda\>
\end{equation}
where the sum is over all partitions $\lambda$  such that $\mu/\lambda$ is
a horizontal $n$-ribbon strip of weight $k$. 

Identifying $\F_q$ with $\Q(q)\otimes Sym$ by setting
$|\lambda\>=s_\lambda$, one can define a linear operator 
$$
\psi^n_q :\ \F_q \longrightarrow \F_q
$$
by specifying the image of the basis $(h_\lambda)$ as
\begin{equation}
\psi_q^n (h_\lambda) =
V_{\lambda_1}V_{\lambda_2}\cdots V_{\lambda_r} |\emptyset\> \ .
\end{equation}
Then, the image $\{\psi_q^n(g_\lambda)\}$ of any basis $\{g_\lambda\}$
of symmetric functions  will be a basis of the space of $U_q(\slchap_n)$-highest
weight vectors in $\F_q$.

%%%%%%%%%%%%%%%%%%%%%%%%%%%%%%%%%%%%%%%%%%%%%%%%%%%%%%%%%%%%%%%%%%%%%%%%%%%%%
%%%%%%%%%%%%%%%%%%%%%%%%%%%%%%%%%%%%%%%%%%%%%%%%%%%%%%%%%%%%%%%%%%%%%%%%%%%%%

\section{An involution of the Fock space}

Let $\I$ be the set of decreasing sequences 
$I=(i_1,i_2,\cdots )$ such that $i_k=-k+1$ for $k$ large enough.
Then $\{u_I\, | \, I \in \I\}$ is the standard basis of 
$\F_q$. 
For $m\ge 0$, denote by $\I_m$ the the subset of $\I$ consisting of those $I$
such that $\sum_k (i_k+k-1) = m$.
Let $I\in \I_m$, and let $u_I$ be the associated basis vector of
$\F_q$. We denote by $\alpha_{n,k}(I)$ the number of pairs
$(r,s)$ with $1\le r < s \le k$ and $r-s \not\equiv 0 \mod n$.

\begin{proposition}\label{prop:3.1}
For $k\ge m$, the $q$-wedge
$$
\overline{u_I} =
(-1)^{{k\choose 2}} q^{\alpha_{n,k}(I)} 
u_{i_k}\wedge u_{i_{k-1}}\wedge \cdots\wedge u_{i_1}\wedge
u_{i_{k+1}}\wedge u_{i_{k+2}} \wedge \cdots
$$
is independent of $k$.
\end{proposition}

Define a semi-linear map $v\mapsto \overline{v}$ in $\F_q$
by
\begin{equation}
\overline{\sum_{I\in\I}\varphi_I(q)u_I} =
\sum_{I\in\I}\varphi_I(q^{-1}) \overline{u_I} \ .
\end{equation}

\begin{theorem}\label{th:3.2}\par
\begin{quote}
{\rm (i)} \quad $v\mapsto \overline{v}$ is an involution of $\F_q$. 

{\rm (ii)} \quad $\overline{f_i v} = f_i \overline{v}$
and $\overline{B_{-k}v}=B_{-k}\overline{v}$, \quad
{\rm (}$v\in\F_q$, $i\in\{0,\ldots,n-1\}$, $k>0${\rm )}.
\end{quote}
\end{theorem}
We note that there is a unique semi-linear map satisfying (ii) and
$\overline{|\emptyset\>}=|\emptyset\>$. This implies that the restriction
of the involution $v\mapsto \overline{v}$ to the subspace $M(\Lambda_0)$
of $\F_q$ coincides with the usual involution in terms of which the
global crystal basis of $M(\Lambda_0)$ is defined.

Let $\mu\vdash m$. Set
$$
\overline{|\mu\>}=\sum_{\lambda\vdash m} a_{\lambda\mu}(q) |\lambda\> \ .
$$
\begin{theorem}\label{th:3.3}
\begin{quote}
{\rm (i)} \quad $a_{\lambda\mu}(q)\in \Z[q,q^{-1}]$. 

{\rm (ii)} \quad $a_{\lambda\mu}(q)=0$ unless $\lambda\unlhd\mu$ and
$\lambda$, $\mu$ have the same $n$-core. 

{\rm (iii)} \quad $a_{\lambda\lambda}(q)=1$. 

{\rm (iv)} \quad $a_{\lambda\mu}(q) = a_{\mu'\lambda'}(q)$.
\end{quote}
\end{theorem}

For $n=2$, the matrices 
${\bf A}_m(q)=[a_{\lambda\,\mu}(q)]_{\lambda,\,\mu \vdash m}$ 
for $m=2,3,4$ are
$$
\left [\begin {array}{cc} 1&0\\q-q^{-1}&1\end {array}
\right ] \qquad
\left [\begin {array}{ccc} 1&0&0\\0&1&0\\q-q^{-1}&0&1
\end {array}\right ]
\qquad
\left [\begin {array}{ccccc} 1&0&0&0&0\\q-q^{-1}&1&0&0&0\\
-1+q^{-2}&q-q^{-1}&1&0&0\\0&q^{2}-1&
q-q^{-1}&1&0\\q^{2}-1&0&-1+q^{-2}&
q-q^{-1}&1\end {array}\right ]
$$
(partitions are ordered in reverse lexicographic order,
e.g., for $m=4$, $(4),(31),(22),(211),(1111)$).

%%%%%%%%%%%%%%%%%%%%%%%%%%%%%%%%%%%%%%%%%%%%%%%%%%%%%%%%%%%%%%%%%%%%%%%%%%%%%%
%%%%%%%%%%%%%%%%%%%%%%%%%%%%%%%%%%%%%%%%%%%%%%%%%%%%%%%%%%%%%%%%%%%%%%%%%%%%%%

\section{Canonical bases}

Let $L$ (resp. $L^-$) be the $\Z[q]$ (resp. $\Z[q^{-1}]$)-lattice in $\F_q$
with basis $\{ |\lambda\>\}$. Using Theorems \ref{th:3.2} and
\ref{th:3.3} one can construct ``IC-bases of $\F_q$'', in the terminology
of Du \cite{Du} (see also \cite{Lu}, 7.10).

\begin{theorem}\label{th:4.1}
There exist bases $\{G(\lambda)\}$, $\{G^-(\lambda)\}$ of $\F_q$
characterized by: 
\begin{quote}
{\rm (i)} \quad $\overline{G(\lambda)}=G(\lambda)$, \quad 
$\overline{G^-(\lambda)}=G^-(\lambda)$,

{\rm (ii)} \quad $G(\lambda)\equiv |\lambda\> \mod qL$,
\quad $G^-(\lambda)\equiv |\lambda\> \mod q^{-1}L^-$.
\end{quote}
\end{theorem}
Set 
$$\displaystyle G(\mu)=\sum_\lambda d_{\lambda\mu}(q)|\lambda\>\,,\quad\quad
\displaystyle G^-(\lambda)=\sum_\mu e_{\lambda\mu}(q) |\mu\>\,.$$ 
Then
$d_{\lambda\mu}(q)\in\Z[q]$, $e_{\lambda\mu}(q)\in\Z[q^{-1}]$, and these
polynomials are nonzero only if 
%$\lambda\unlhd\mu$ and 
$\lambda,\mu$ have
the same $n$-core. 
Moreover, $d_{\lambda\lambda}(q)=e_{\lambda\lambda}(q)=1$,
and $d_{\lambda\mu}(q)=0$ unless $\lambda\unlhd\mu$,
$e_{\lambda\mu}(q)=0$ unless $\mu\unlhd\lambda$.
Also, $\{G(\lambda)\,|\, \lambda \ n\mbox{\rm -regular}\}$ coincides
with the lower crystal basis of the basic representation $M(\Lambda_0)$
of $U_q(\slchap_n)$.

Let $\{G^\dagger(\lambda)\}$ be the adjoint basis of $\{G(\lambda)\}$.
It follows from Theorem \ref{th:3.3} (iv) that
$$
G^\dagger (\lambda)' = G^-(\lambda')
$$
where $v\mapsto v'$ denotes the semi-linear involution of $\F_q$
defined by $|\lambda\>'=|\lambda'\>$. 

We set $$\displaystyle G^\dagger(\lambda)=\sum_\mu c_{\lambda\mu}(q) |\mu\>\,,$$
so that
$c_{\lambda\mu}(q)=e_{\lambda'\mu'}(q^{-1})$ and 
$[c_{\lambda\mu}(q)] =[d_{\lambda\mu}(q)]^{-1}$.

For $n=2$ and $m\le 6$, the matrices 
$\D_m(q)=[d_{\lambda\mu}(q)]_{\lambda,\,\mu \vdash m}$
are

$$
\begin {array}{ccc}  2 &1&0\\ 1 1 &q&1\end {array}  \qquad\qquad
\begin {array}{cccc}  3 &1&0&0\\ 2 1 &0&1&0\\ 1 1 1 &q&0&1
\end {array}  
$$

\vskip 10mm

$$
\begin {array}{cccccc}  4 &1&0&0&0&0\\ 3 1 &q&1&0&0&0\\ 2 2 &0&
q&1&0&0\\ 2 1 1 &q&q^{2}&q&1&0\\ 1 1 1 1 &q^{2}&0&0&q&1\end {array}
 \qquad\quad 
  \begin {array}{cccccccc}  5 &1&0&0&0&0&0&0\\ 4 1 &0&1&0&0&0&0&0
\\ 3 2 &0&0&1&0&0&0&0\\ 3 1 1 &q&0&q&1&0&0&0\\ 2 2 1 &0&0&q^{2}&q&1&0&0
\\ 2 1 1 1 &0&q&0&0&0&1&0\\ 1 1 1 1 1 &q^{2}&0&0&q&0&0&1\end {array}
$$  

\vskip 10mm
$$
\begin {array}{cccccccccccc}  6 &1&0&0&0&0&0&0&0&0&0&0\\ 5 1 &q
&1&0&0&0&0&0&0&0&0&0\\ 4 2 &0&q&1&0&0&0&0&0&0&0&0\\ 4 1 1 &q&q^{2}&q&1
&0&0&0&0&0&0&0\\ 3 3 &0&0&q&0&1&0&0&0&0&0&0\\ 3 2 1 &0&0&0&0&0&1&0&0&0
&0&0\\ 3 1 1 1 &q^{2}&q&q^{2}&q&q&0&1&0&0&0&0\\ 2 2 2 &0&0&q^{2}&q&q&0
&0&1&0&0&0\\ 2 2 1 1 &0&q^{2}&q^{3}&q^{2}&q^{2}&0&q&q&1&0&0\\ 2 1 1 1 
1 &q^{2}&q^{3}&0&q&0&0&q^{2}&0&q&1&0\\ 1 1 1 1 1 1 &q^{3}&0&0&q^{2}&0&0
&0&0&0&q&1\end {array}  
$$

%%%%%%%%%%%%%%%%%%%%%%%%%%%%%%%%%%%%%%%%%%%%%%%%%%%%%%%%%%%%%%%%%%%%%%%%%%
%%%%%%%%%%%%%%%%%%%%%%%%%%%%%%%%%%%%%%%%%%%%%%%%%%%%%%%%%%%%%%%%%%%%%%%%%%

\section{A $q$-analogue of Steinberg's tensor product theorem}

Let $\alpha$ be a partition of $r$. Define an element $S_\alpha$ of $\U$ by
\begin{equation}
S_\alpha = \sum_{\beta\vdash r}{\chi^\alpha_\beta\over z_\beta}
B_{-\beta_1}B_{-\beta_2}\cdots B_{-\beta_k}
\end{equation}
where $\beta=(\beta_1,\ldots,\beta_k)=(1^{m_1}2^{m_2}\ldots r^{m_r})$,
$z_\beta= 1^{m_1} m_1! \cdots r^{m_r} m_r!$ and $\chi^\alpha_\beta$
is the value of the irreducible character $\chi^\alpha$ of the symmetric group
on a permutation of cycle type $\beta$. 

Writing $V_\mu = V_{\mu_1}V_{\mu_2}\cdots V_{\mu_r}$, one has
$$
S_\alpha = \sum_{\mu\vdash r} \kappa_{\alpha \mu} V_\mu\,,
$$
where the $\kappa_{\alpha \mu}$ denote the entries of the inverse
Kostka matrix. Hence, using (\ref{eq8}), one can describe the
action of $S_\alpha$ on Young diagrams.

The operator $S_\alpha$ is a $q$-analogue of the multiplication by
the plethysm $\psi^n(s_\alpha)$ in the ring of symmetric functions.

Let $\lambda$ be a partition such that $\lambda'$ is $n$-singular. We can write
$\lambda=\mu+n\alpha$ where $n\alpha=(n\alpha_1,n\alpha_2,\ldots)$ and
$\mu'$ is $n$-regular.

\begin{theorem}\label{th:5.1}
\quad $G^-(\lambda) = S_\alpha\left( G^-(\mu)\right)$.
\end{theorem}
This reduces the computation of $\{G^-(\lambda)\}$ to that of the
subfamily $\{G^-(\mu)\,|\,\mu\ n\mbox{\rm -regular}\}$.

Let $\D_m$ denote the decomposition matrix of the $v$-Schur algebra ${\cal S}_m$
over a field of characteristic $0$, for $v$ a primitive $n$-th root
of unity.
We use the notational convention of James \cite{Ja}, that is,
the rows and columns of $\D_m$ are indexed in such a way that
$\D_m$ is the matrix $\Delta_m$ of \cite{Ja} for big $p$.

\begin{conjecture}\label{conj:5.2}
The matrix $\D_m(1)$ is equal to the decomposition matrix $\D_m$. 
\end{conjecture}
This conjecture, which generalizes Conjecture 6.9 of \cite{LLT2}, is
already verified to a large extent. Indeed, on the one hand Ariki \cite{A}
and Grojnowski have verified independently
our previous conjecture, which means that the $d_{\lambda\mu}(1)$ for
$\mu$ $n$-regular are equal to the corresponding decomposition numbers of
${\cal S}_m$. On the other hand, it follows from results of James \cite{Ja}
that the entries of the inverse matrices ${\bf D}_m^{-1}$ satisfy the same
properties as those deduced from Theorem~\ref{th:5.1} for the
coefficients $c_{\lambda\mu}(1)$. Since we have verified, using the tables
of \cite{Ja}, that
${\bf D}_m^{-1} = [c_{\lambda\mu}(1)]$ for $m\le 10$ (any~$n$),
we deduce that an infinite number of $c_{\lambda\mu}(1)$ coincide with the
corresponding entries of ${\bf D}_m^{-1}$.

The following refined conjecture has also been checked for small $m$. It
generalizes the conjecture of Section 9 of \cite{LLT2}, due to Rouquier.
Let $(W(\lambda)^i)$ be the Jantzen filtration of the Weyl module
$W(\lambda)$ for the $v$-Schur algebra ${\cal S}_m$, and let $L(\mu)$ be the
irreducible module corresponding to $\mu$ \cite{JM2}.
\begin{conjecture} Let $\lambda$, $\mu$ be partitions of $m$. Then,
$$
d_{\lambda'\mu'}(q)=\sum_{i\ge 0} [W(\lambda)^i/W(\lambda)^{i+1}:L(\mu)]\,q^i\,.
$$
\end{conjecture}

Finally, we note the following combinatorial description of some polynomials
$e_{\lambda\mu}(q)$ in the case $n=2$, which proves that they are equal
(up to sign and the replacement of $q$ by $q^{-2}$) 
to the $q$-analogues of Littlewood-Richardson coefficients introduced
in~\cite{CL}.

\begin{theorem}
Let $n=2$. One has
$$
e_{2\lambda,\mu}(q) =\varepsilon_2(\mu)\sum_{T\in{\rm Yam}_2(\mu,\lambda)}
q^{-2{\rm spin}\,(T)}
$$
where ${\rm Yam}_2(\mu,\lambda)$ denotes the set of Yamanouchi domino tableaux
of shape $\mu$ and weight $\lambda$, and $\varepsilon_2(\mu)$ is the 2-sign of 
$\mu$.
\end{theorem}

%%%%%%%%%%%%%%%%%%%%%%%%%%%%%%%%%%%%%%%%%%%%%%%%%%%%%%%%%%%%%%%%%%%%%%%%%%
%%%%%%%%%%%%%%%%%%%%%%%%%%%%%%%%%%%%%%%%%%%%%%%%%%%%%%%%%%%%%%%%%%%%%%%%%%%


\footnotesize
\begin{thebibliography}{ABCD}

\bibitem{A} {\sc S. Ariki}, {\it On the decomposition numbers of the
Hecke algebra of $G(m,1,n)$}, Preprint 1996.

\bibitem{CL}{\sc C. Carr\'e} and {\sc B. Leclerc},
{\it Splitting the square of a Schur function into its symmetric and 
antisymmetric parts}, J. Alg. Comb. {\bf 4} (1995), 201-231.

\bibitem{DJKM} {\sc E. Date, M. Jimbo, M. Kashiwara} and {\sc  T. Miwa},
{\it Transformation groups for soliton equations},
in ``Non-linear Integrable Systems -- Classical Theory
and Quantum Theory", Proceedings of RIMS Symposium,
World Scientific 1983, 39-119.

\bibitem{DJKMO} {\sc E. Date, M. Jimbo, A. Kuniba, T. Miwa} and {\sc M. Okado},
{\it Paths, Maya diagrams, and representations of $\slchap(r,C)$},
Adv. Stud. Pure Math. {\bf 19} (1989), 149-191.

\bibitem{DF} {\sc J. Ding} and {\sc I.B. Frenkel}, {\it Isomorphism of
two realizations of quantum affine $U_q(\glchap_n)$},
Commun. Math. Phys. {\bf 156} (1993), 277-300.

\bibitem{Du} {\sc J. Du}, {\it IC bases and quantum linear spaces}, 
Proc. Symp. Pure Math. {\bf 56} (1994), part 2, 135-148. 

\bibitem{GRV} {\sc V. Ginzburg, N. Reshetikhin} and {\sc E. Vasserot},
{\it Quantum groups and flag varieties}, Contemp. Math. {\bf 175}
(1994), 101-130.

\bibitem{Ha} {\sc T. Hayashi}, {\it $q$-analogues of Clifford and Weyl algebras -
spinor and oscillator representations of quantum enveloping algebras},
Commun. Math. Phys., {\bf 127} (1990), 129-144.

\bibitem{Ja} {\sc G. James}, {\it The decomposition matrices 
of $GL_n(q)$ for $n\le 10$},
Proc. London Math. Soc., {\bf 60} (1990), 225-265.

\bibitem{JM2} {\sc G. James} and {\sc  A. Mathas}, {\it A $q$-analogue of 
the Jantzen-Schaper theorem}, Preprint 1995.

\bibitem{KMS}{\sc M. Kashiwara, T. Miwa} and {\sc E. Stern}, {\it Decomposition
of $q$-deformed Fock spaces}, preprint RIMS-1035 (1995).

\bibitem{LLT1} {\sc A. Lascoux, B. Leclerc} and {\sc  J.-Y. Thibon}, 
{\it Une conjecture
pour le calcul des matrices de d\'ecomposition des alg\`ebres de Hecke de
type $A$ aux racines de l'unit\'e}, C. R. Acad. Sci. Paris {\bf 321} (1995),
511-516.

\bibitem{LLT2} {\sc A. Lascoux, B. Leclerc} and {\sc  J.-Y. Thibon}, {\it Hecke algebras
at roots of unity and crystal bases of quantum affine algebras},
Commun. Math. Phys. (to appear).

\bibitem{LLT3} {\sc A. Lascoux, B. Leclerc} and {\sc  J.-Y. Thibon}, {\it Ribbon
tableaux, Hall-Littlewood functions, quantum affine algebras and unipotent
varieties}, preprint q-alg$/$9512031.

\bibitem{Lu}{\sc G. Lusztig}, {\it Canonical bases arising from quantized
enveloping algebras}, J. Amer. Math. Soc. {\bf 3} (1990), 447-498.

\bibitem{Mcd}{\sc I.G. Macdonald}, {\it Symmetric functions and Hall polynomials},
2nd edition, Oxford 1995.

\bibitem{MM}{\sc K.C. Misra } and {\sc T. Miwa}, {\it Crystal base for the basic
representation of $U_q(\slchap_n)$}, Commun. Math. Phys. {\bf 134} (1990),
79-88.

\bibitem{St}{\sc E. Stern}, {\it Semi-infinite wedges and vertex operators},
Internat. Math. Research Notices {\bf 4} (1995), 201-220.

\end{thebibliography}
\end{document}